\documentclass[12pt]{amsart}
\usepackage{amsthm,amssymb}

\newtheorem{thm}{Theorem}[section]
\newtheorem{lm}[thm]{Lemma}
\newtheorem{cor}[thm]{Corollary}
\newtheorem{pro}[thm]{Proposition}
\newtheorem{pro-def}[thm]{Proposition and Definition}

\theoremstyle{definition}

\input epsf

\begin{document}

\title[]{On the bridge number of knot diagrams with minimal crossings}
\author[]{Jae-Wook Chung and Xiao-Song Lin}
\address{Department of Mathematics, University of California, Riverside, CA 92521}
\email{xl@math.ucr.edu, jwchung@math.ucr.edu}
\thanks{The second named author is supported in part by NSF}

\begin{abstract}{Given a diagram $D$ of a knot $K$, we consider the
number $c(D)$ of crossings and the number $b(D)$ of overpasses of
$D$. We show that, if $D$ is a diagram of a nontrivial knot $K$
whose number $c(D)$ of crossings is minimal, then $1+\sqrt{1+c(D)}
\leq b(D)\leq c(D)$. These inequalities are shape in the sense
that the upper bound of $b(D)$ is achieved by alternating knots
and the lower bound of $b(D)$ is achieved by torus knots. The
second inequality becomes an equality only when the knot is an
alternating knot. We prove that the first inequality becomes an
equality only when the knot is a torus knot.}
\end{abstract}

\maketitle

\section{Introduction}

A diagram $D$ of a knot $K$ is a nice representative of the knot
type $[K]$ which is the isotopy class of $K$. It can be obtained
from a regular planar projection (or simply, regular projection)
$P$ of $K$ as follows. Let us take a sufficiently small
neighborhood of each double point of $P$ so that the intersection
of the neighborhood and $P$ is of the `shape X' on the plane. Then
modify the interior of each neighborhood so that we get a knot $D$
which is isotopic to $K$ and regularly projected to $P$. In this
sense, a knot diagram can be `almost planar', i.e., it lies in the
plane except for a sufficiently small neighborhood of each double
point of the regular planar projection.

Tabulation of knots is usually done by listing knots using their
diagrams of minimal numbers of crossings. In general, such a
minimal knot diagram is not unique within its isotopy class. Can
this nonuniqueness be fixed by considering some other quantities
associated with knot diagrams? In this paper, we study the number
of overpasses of a minimal knot diagram. An immediate concern in
this investigation is the relationship between the number of
crossings and the number of overpasses of a knot diagram. As it
turns out, the number of crossings can be estimated from below and
above by that of overpasses if the number of crossings is minimal
among all diagrams of the same knot type (Theorem 2.9).

As an interesting consequence of this estimation of the minimal
crossing number by the number of overpasses, we may put two nice
and relatively well understood classes of knots, that of
alternating knots and that of $(k-1,\pm k)$-torus knots, on the
opposite extremes of a measurement of nonalternatingness of knot
diagrams (Theorem 2.10). Under this measurement, we prove that the
`most nonalternating' knot diagrams are exactly the standard
minimal $(k-1,\pm k)$-torus knot diagrams (Theorem 3.1).

Notice that the minimal number of overpasses of a knot was the
classical {\it bridge number} of a knot, first studied by Schubert
in \cite{SCH}, where the effect of various operations on knots
(satellite, cabling, connected sum) on this number was
investigated.

Throughout this paper, all knots are oriented and lie in the
3-dimensional sphere $S^3$. Also, all knots are tame, i.e., they
are isotopic to polygonal or smooth knots. Hence, every knot has a
diagram whose number of crossings is finite. For convenience, we
distinguish knot diagrams and knots. From now on, a knot means an
isotopy knot type $[D]$ for some knot diagram $D$. We will refer
the reader to textbooks of knot theory (such as \cite{A} or
\cite{BZ}) for standard terminologies. Also, we assume the
\textit{well ordering property} of the set $N \cup \{0\}$ of
nonnegative integers to guarantee existence of the smallest
element of any nonempty subset of $N \cup \{0\}$.

We would like to thank Professor Marek Chrobak of the Department
of Compute Science and Engineering at UC Riverside for asking
stimulating questions which led our inequality about the numbers
$c(D)$ and $b(D)$ to its current optimal form.

\section{Minimal crossings and bridge number of knot diagrams}

The number of crossings of a knot diagram $D$ is denoted by
$c(D)$. For each knot $K$, we denote $min\{c(D) \mid D$ \textit{is
a diagram of} $K\}$ by $c(K)$. Note that we may assume a knot
diagram $D$ lies in the plane by indicating `overcrossings' and
`undercrossings'. A `crossing', in normal sense, of a knot diagram
$D$ means a `double point' of the regular projection of $D$.
Hence, $c(D)$ is the number of all double points of the regular
projection of $D$. If $x$ is a crossing of a knot diagram $D$, we
denote the overcrossing and the undercrossing of $D$ projected to
$x$ by $x_+$ and $x_-$, respectively, where $x_+$ is above $x$ and
$x_-$ is below $x$. Alternatively, we may regard a crossing of $D$
as a pair of two points, overcrossing and undercrossing, in $D$
which are projected to the same double point. In this case, a
crossing is considered as the preimage of a double point under the
projection map.

\begin{pro} Let $D$ be a knot diagram with $c(D) \geq 1$. Then
there is a unique positive integer $k$ such that there is a finite
sequence $$s_1,f_1,s_2,f_2,...,s_k,f_k$$ of $2k$ points of $D$,
each of which is neither an overcrossing nor an undercrossing of
$D$, such that
$$[s_1,f_1],[s_2,f_2],...,[s_{k-1},f_{k-1}],[s_k,f_k]$$ and
$$[f_1,s_2],[f_2,s_3],...,[f_{k-1},s_k],[f_k,s_1]$$ are the
{\em overpasses} and the {\em underpasses} of $D$ with respect to
$s_1,f_1,s_2,f_2,...,s_k,f_k$, respectively, where $[s_i, f_i]$,
for each $i \in \{1,...,k\}$, is the closed arc of $D$ from $s_i$
to $f_i$ which contains at least one overcrossing but has no
undercrossing; similarly, $[f_i,s_{i+1}]$, for each $i \in
\{1,...,k\}$, with the subscript counted modulo $k$, is the closed
arc of $D$ from $f_i$ to $s_{i+1}$ which contains at least one
undercrossing but has no overcrossing.
\end{pro}

\begin{proof} We may assume $c(D) \geq 2$. If $c(D)=1$,
then $D$ is of the `shape 8', and hence, $k=1$.

Existence: We will construct a finite sequence as stated above.
The following procedure is one way to get it. Let us start at a
point $\ast$ on $D$ slightly before an undercrossing and go along
with $D$ until arriving at the first overcrossing $a_1$ from
$\ast$. Let $a_2$ be the undercrossing just before $a_1$. Take a
point $s_1$ of $D$ between $a_2$ and $a_1$. Then go along with $D$
from $s_1$ until arriving at the first undercrossing $a_3$ from
$s_1$, and let $a_4$ be the overcrossing just before $a_3$. Take a
point $f_1$ of $D$ between $a_4$ and $a_3$. We may now repeat this
procedure to take the other points $s_2,f_2,s_3,f_3,...,s_k,f_k$
until there is no overcrossing between $f_k$ and $\ast$. Since the
number $c(D)$ of crossings of $D$ is finite, $k$ must be finite.
By its construction, the sequence of points
$s_1,f_1,s_2,f_2,...,s_k,f_k$ of $D$ is a one as stated in the
proposition.

Uniqueness: Suppose that $l$ is a positive integer and
$s'_1,f'_1,s'_2,f'_2,...,s'_l,f'_l$ is a sequence of $D$ as stated
in the proposition. Let
$$a_2,a_1,a_4,a_3,...,a_{2(2k-1)},a_{2(2k-1)-1},a_{2(2k)},a_{2(2k)-1}$$
be the sequence of points on $D$ as described above. Then $s'_1$
must be contained in an arc of $D$ which contains $s_i$ for some
$i \in \{1,...,k\}$ but no overcrossing and no undercrossing of
$D$. In other words, $s'_1$ must be between $a_{2(2i-1)}$ and
$a_{2(2i-1)-1}$ for some $i \in \{1,...,k\}$. Hence, $f'_1$ must
be between $a_{2(2i)}$ and $a_{2(2i)-1}$. Keep going on like this.
If $i=1$, then $f'_l$ must be between $a_{2(2k)}$ and
$a_{2(2k)-1}$, and if $1<i \leq k$, then $f'_l$ must be between
$a_{2(2(i-1))}$ and $a_{2(2(i-1))-1}$. Hence, we have $l=k$. This
proves the uniqueness of $k$.
\end{proof}

Such a sequence $s_1,f_1,s_2,f_2,...,s_k,f_k$ as in Proposition
2.1 is said to be an \textit{over-underpass sequence} of $D$.
Since any over-underpass sequence of $D$ consists of $2k$ points,
the number of overpasses (or underpasses) with respect to any
over-underpass sequence of $D$ is $k$. Hence, we define the number
of overpasses (or underpasses) of the knot diagram $D$ as $k$, and
denote it by $b(D)$, sometimes called the length of over-underpass
sequence. If a knot diagram $D$ has no crossing, we define $b(D)$
as 0. Also, for each knot $K$, we denote $min\{b(D) \mid D$
\textit{is a diagram of} $K\}$ by $b(K)$. This number $b(K)$ is
called the \textit{bridge number} of $K$.

Notice that (1) $c(D)$ and $b(D)$ are plane isotopy invariants of
knot diagrams, i.e., if $D_1$ and $D_2$ are plane isotopic knot
diagrams, then $c(D_1)=c(D_2)$ and $b(D_1)=b(D_2)$; (2) $c(K)$ and
$b(K)$ are isotopy invariants of knots, i.e., if $K_1$ and $K_2$
are isotopic knots, then $c(K_1)=c(K_2)$ and $b(K_1)=b(K_2)$.

\begin{cor} If $b(D)\leq 1$, then $D$ is a diagram of a
trivial knot. Therefore, a knot $K$ is trivial if and only if $K$
has a diagram $D$ with $b(D)\leq 1$.
\end{cor}

Obviously, for any positive integer $k$, there is a diagram of a
trivial knot whose number of overpasses is greater than $k$. On
the other hand, given a knot diagram $D$ with at least one
crossing, we can add crossings to $D$ as many as we want without
changing the knot type and the number of overpasses of $D$. For
example, take a sufficiently small arc of $D$ from $s_1$ to a
point between $s_1$ and the first overcrossing of $D$ from $s_1$,
and twist it alternatively so that the number of overpasses of $D$
is not changed. Or, we may modify the interior of a sufficiently
small neighborhood of a crossing of $D$ to achieve the goal of
increasing $c(D)$ arbitrarily while keeping $b(D)$ and the knot
type of $D$ fixed. Hence, we have the following corollary.

\begin{cor} If $D$ is a diagram of a knot $K$ such that
$c(D)\geq 1$, then for every positive integer $n$, there is a
diagram $D'$ of $K$ such that $b(D')=b(D)$ and $c(D') \geq
c(D)+n$.
\end{cor}

Remark that the number of crossings of a knot diagram with a
minimal number of overpasses can be arbitrarily large.

\begin{lm} $b(D)\leq c(D)$ for any knot diagram $D$. The
equality holds if and only if $D$ is an alternating knot diagram.
Furthermore, $b(K)\leq c(K)$ for any knot $K$.
\end{lm}

\begin{proof} Suppose that $b(D)=k \geq 1$ and
$s_1,f_1,s_2,f_2,...,s_k,f_k$ is an over-underpass sequence of
$D$. Let $m_i$ be the number of overcrossings of $D$ on
$[s_i,f_i]$ for each $i \in \{1,...,k\}$, $n_i$ the number of
undercrossings of $D$ on $[f_i,s_{i+1}]$ for each $i \in
\{1,...,k-1\}$, and $n_k$ the number of undercrossings of $D$ on
$[f_k,s_1]$. Then
$$2c(D)=\sum_{i=1}^k m_i + \sum_{i=1}^k n_i \geq
\sum_{i=1}^k 1 + \sum_{i=1}^k 1=2b(D).$$ $D$ is an alternating
knot diagram if and only if every overpass and underpass has
exactly one overcrossing and undercrossing, respectively, i.e.,
$m_i=n_i=1$ for each $i \in \{1,...,k\}$.

Also, it follows immediately that $b(K)\leq c(K)$ for any knot
$K$.
\end{proof}

The following lemma is the first key to prove the first main
theorem of this paper (Theorem 2.9). By this lemma, if an overpass
of a knot diagram $D$ with respect to an over-underpass sequence
of $D$ crosses an underpass more than once, then the number of
crossings of $D$ is not minimal any more.

\begin{lm} If $D$ is a minimal diagram of a knot $K$ with respect to
crossings, i.e., $c(D)=c(K)$, $b(D)$ is a nonnegative integer $k$,
and $s_1,f_1,s_2,f_2,...,s_k,f_k$ is an over-underpass sequence of
$D$, then every overpass of $D$ with respect to
$s_1,f_1,s_2,f_2,...,s_k,f_k$ crosses each underpass at most once.
\end{lm}

\begin{proof} We may assume $K$ is nontrivial. If $K$ is
trivial, then $c(D)=c(K)=0$, hence, $b(D)=0$. Let $D$ be a diagram
of a knot $K$ such that $c(D)=c(K)$. Suppose that $b(D)=k \geq 1$
and $s_1,f_1,s_2,f_2,...,s_k,f_k$ is an over-underpass sequence of
$D$. Then, by Lemma 2.4, $c(D) \geq k \geq 1$. Let $o_i=[s_i,f_i]$
and $u_i=[f_i,s_{i+1}]$ for each $i \in \{1,...,k-1\}$, and let
$o_k=[s_k,f_k]$ and $u_k=[f_k,s_1]$. Suppose that the number of
crossings between $o_i$ and $u_j$ is at least 2 for some $i,j \in
\{1,...,k\}$. Let $x_1$ and $x_2$ be two crossings of $D$ such
that $x_{1+}$ and $x_{2+}$ are the first and the second
overcrossings from $s_i$ among all overcrossings between $o_i$ and
$u_j$, respectively, and let $y_1$ and $y_2$ be two crossings of
$D$ such that $y_{1-}$ and $y_{2-}$ are the first and the second
undercrossings from $f_j$ in $\{x_{1-},x_{2-}\}$, respectively.
Then either $y_1=x_1$, $y_2=x_2$ or $y_1=x_2$, $y_2=x_1$. There
are two different diagrams for each case. One of them for the case
of $y_1=x_1$, $y_2=x_2$ is shown in Figure 1. It should be not
hard for the reader to figure out the other diagrams.

Let $r$ be the number of overcrossings on the arc
$\overrightarrow{x_{1+}x_{2+}}^{o_i}$ of $o_i$ from $x_{1+}$ to
$x_{2+}$, and let $s$ be the number of undercrossings on the arc
$\underrightarrow{y_{1-}y_{2-}}_{u_j}$ of $u_j$ from $y_{1-}$ to
$y_{2-}$. To prove this lemma, for each of the cases $r \leq s$
and $r>s$, we will construct a diagram $D'$ of $K$ by modifying an
arc of $o_i$ or $u_j$ such that $c(D')<c(D)$. Temporarily, we
regard the knot diagram $D$ as its regular projection. Hence, $D$
lies in the plane, and all overcrossings and all undercrossings of
$D$ are the double points of $D$ as the regular projection.
Consider the $r$ double points $d_1,d_2,...,d_r$ on the arc
$\overrightarrow{x_1x_2}^{o_i}$ of $o_i$ from $x_1$ to $x_2$ such
that $d_1<d_2<...<d_r$ with respect to the order from $s_i$ and
the $s$ double points $e_1,e_2,...,e_s$ on the arc
$\underrightarrow{y_1y_2}_{u_j}$ of $u_j$ from $y_1$ to $y_2$ such
that $e_1<e_2<...<e_s$ with respect to the order from $f_j$. Then
$x_1=d_1$, $x_2=d_r$ and $y_1=e_1$, $y_2=e_s$. Notice that we can
take a sufficiently small positive real number $\epsilon$, an
$\epsilon$-neighborhood $U_{o_i,\epsilon}$ of
$\overrightarrow{x_1x_2}^{o_i}$, and an $\epsilon$-neighborhood
$V_{u_j,\epsilon}$ of $\underrightarrow{y_1y_2}_{u_j}$ such that
$\overline{U_{o_i,\epsilon}} \cap
\{s_1,f_1,s_2,f_2,...,s_k,f_k\}=\varnothing$,
$\overline{V_{u_j,\epsilon}} \cap
\{s_1,f_1,s_2,f_2,...,s_k,f_k\}=\varnothing$, the set of all
double points of $D$ contained in $\overline{U_{o_i,\epsilon}}$ is
$\{d_1,d_2,...,d_r\}$, the set of all double points of $D$
contained in $\overline{V_{u_j,\epsilon}}$ is
$\{e_1,e_2,...,e_s\}$, and, for every positive real number
$\epsilon' \leq \epsilon$, $|Bd(U_{o_i,\epsilon'}) \cap D|=2(r+1)$
and $|Bd(V_{u_j,\epsilon'}) \cap D|=2(s+1)$, where
$Bd(U_{o_i,\epsilon'})$ and $Bd(V_{u_j,\epsilon'})$ are the
boundaries of $U_{o_i,\epsilon'}$ and $V_{u_j,\epsilon'}$,
respectively. Hence, we may assume that
$\overline{U_{o_i,\epsilon}} \cap D=\overrightarrow{a_1a_2}^{o_i}
\cup (l_1 \cup l_2 \cup ... \cup l_r)$ and
$\overline{V_{u_j,\epsilon}} \cap D=\underrightarrow{b_1b_2}_{u_j}
\cup (l^1 \cup l^2 \cup ... \cup l^s)$, where $a_1$ and $a_2$ are
the first and the second points from $s_i$ in
$Bd(U_{o_i,\epsilon}) \cap o_i$, respectively, $b_1$ and $b_2$ are
the first and the second points from $f_j$ in
$Bd(V_{u_j,\epsilon}) \cap u_j$, respectively, $l_t$ is an arc of
the underpass passing through $d_t$ whose endpoints are on
$Bd(U_{o_i,\epsilon})$ for each $t \in \{1,...,r\}$, and $l^t$ is
an arc of the overpass passing through $e_t$ whose endpoints are
on $Bd(V_{u_j,\epsilon})$ for each $t \in \{1,...,s\}$. Remark
that $\{l_1,...,l_r\}$ is pairwise disjoint, $l_t \cap
o_i=\{d_t\}$ and $|l_t \cap Bd(U_{o_i,\epsilon})|=2$ for each $t
\in \{1,...,r\}$, and $|Bd(U_{o_i,\epsilon}) \cap o_i|=2$;
similarly, $\{l^1,...,l^s\}$ is pairwise disjoint, $l^t \cap
u_j=\{e_t\}$ and $|l^t \cap Bd(V_{u_j,\epsilon})|=2$ for each $t
\in \{1,...,s\}$, and $|Bd(V_{u_j,\epsilon}) \cap u_j|=2$.

Let $p_1$, $p_2$, $p_3$, $p_4$ be the first, the second, the
third, the fourth points from $f_j$ in $Bd(U_{o_i,\epsilon}) \cap
(l_1 \cup l_r)$, respectively, and let $\alpha$ be a point on
$\underrightarrow{p_1p_2}_{u_j}$ between $p_1$ and $y_1$. Draw an
arc $A$ in $U_{o_i,\epsilon}$ starting at $\alpha$ so that $A$
intersects $l_t$ transversely only once for each $t \in
\{1,...,r\}$ but $A$ does not intersect $o_i$. Then we take
$\beta$ as the intersecting point of $A$ and
$\underrightarrow{p_3p_4}_{u_j}$ and denote the arc of $A$ from
$\alpha$ to $\beta$ by $\overline{\alpha\beta}$. Similarly,

let $q_1$, $q_2$, $q_3$, $q_4$ be the first, the second, the
third, the fourth points from $s_i$ in $Bd(V_{u_j,\epsilon}) \cap
(l^1 \cup l^s)$, respectively, and let $\gamma$ be a point on
$\overrightarrow{q_1q_2}^{o_i}$ between $q_1$ and $x_1$. Draw an
arc $B$ in $V_{u_j,\epsilon}$ starting at $\gamma$ so that $B$
intersects $l^t$ transversely only once for each $t \in
\{1,...,s\}$ but $B$ does not intersect $u_j$. Then we take
$\delta$ as the intersecting point of $B$ and
$\overrightarrow{q_3q_4}^{o_i}$ and denote the arc of $B$ from
$\gamma$ to $\delta$ by $\overline{\gamma\delta}$.

From now on, $D$ is the diagram of $K$ again, i.e., $D$ is not a
regular projection of a knot but a knot diagram. Let
$\underrightarrow{\alpha\beta}$ be a regularly projecting arc in
$R^3$ whose regular projection is $\overline{\alpha\beta}$ such
that $\underrightarrow{\alpha\beta}$ has no overcrossing, and let
$\overrightarrow{\gamma\delta}$ be a regularly projecting arc in
$R^3$ whose regular projection is $\overline{\gamma\delta}$ such
that $\overrightarrow{\gamma\delta}$ has no undercrossing. Then

\bigskip
\centerline{\epsfxsize=4in\epsfbox{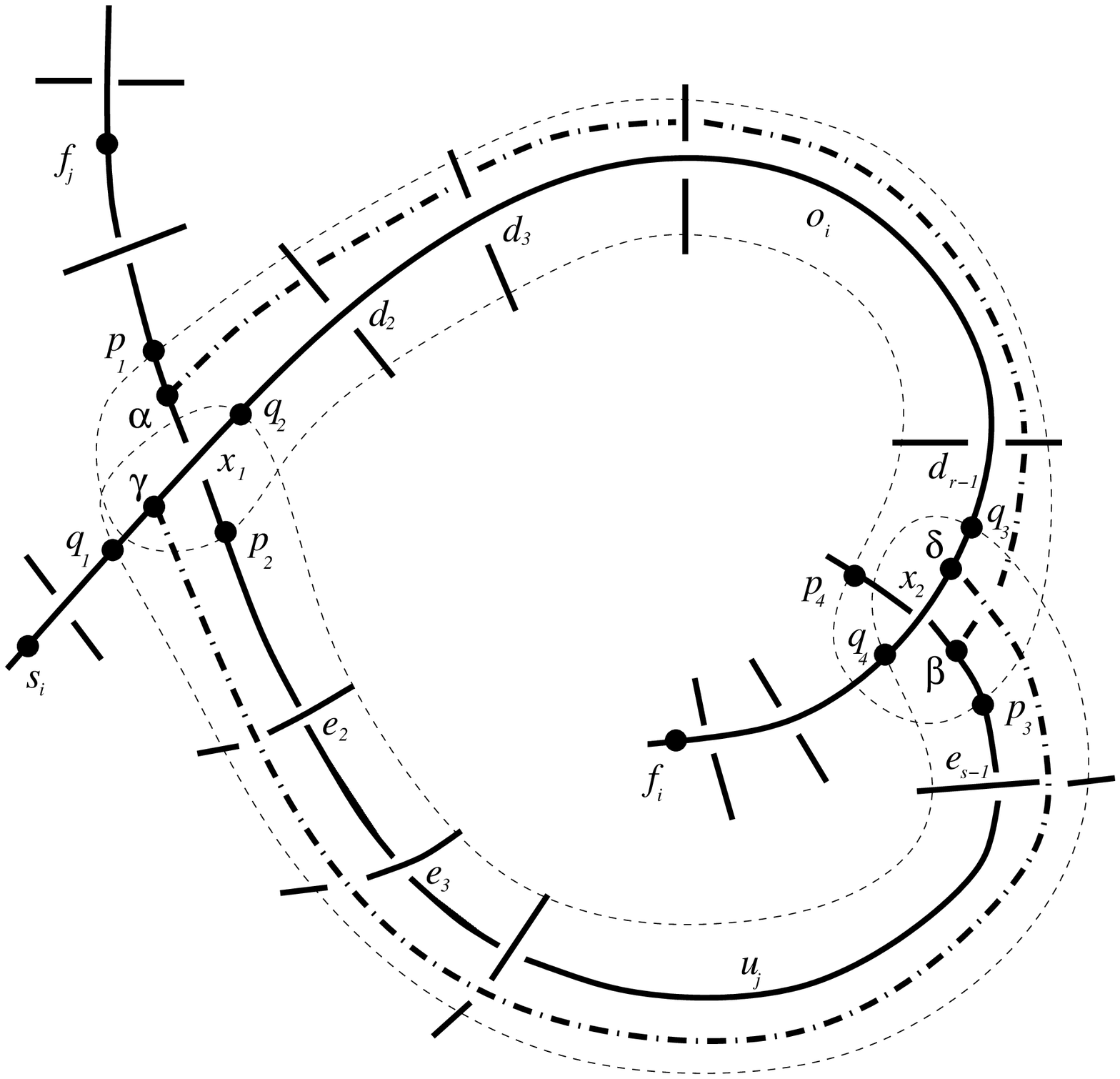}}
\medskip
\centerline{\small Figure 1. Proof of Lemma 2.5: Case 1.(1) and
Case 2.(2) when $y_1=x_1$, $y_2=x_2$.}
\bigskip

\bigskip
\bigskip

Case 1. $r \leq s$: Let $D'=(D-\underrightarrow{\alpha\beta
}_{u_j}) \cup \underrightarrow{\alpha\beta}$. Then $D'$ is a
diagram of $K$.

(1) If $\underrightarrow{p_1p_2}_{u_j}$ and
$\underrightarrow{p_3p_4}_{u_j}$ make the same sign with $o_i$,
then

$c(D')=c(D)-(s-1)+(r-2)$, hence, $c(D)=c(D')+(s-r)+1>c(D')$;

(2) If $\underrightarrow{p_1p_2}_{u_j}$ and
$\underrightarrow{p_3p_4}_{u_j}$ make opposite signs with $o_i$,
then

$c(D')=c(D)-s+(r-2)$, hence, $c(D)=c(D')+(s-r)+2>c(D')$.

Case 2. $r>s$: Let $D'=(D-\overrightarrow{\gamma\delta}^{o_i})
\cup \overrightarrow{\gamma\delta}$. Then $D'$ is a diagram of
$K$.

(1) If $\overrightarrow{q_1q_2}^{o_i}$ and
$\overrightarrow{q_3q_4}^{o_i}$ make the same sign with $u_j$,
then

$c(D')=c(D)-(r-1)+(s-2)$, hence, $c(D)=c(D')+(r-s)+1>c(D')$;

(2) If $\overrightarrow{q_1q_2}^{o_i}$ and
$\overrightarrow{q_3q_4}^{o_i}$ make opposite signs with $u_j$,
then

$c(D')=c(D)-r+(s-2)$, hence, $c(D)=c(D')+(r-s)+2>c(D')$.

Therefore, we have a diagram $D'$ of $K$ such that $c(D')<c(D)$.
This is a contradiction to $c(D)=c(K)$.
\end{proof}

Now, we prove the second key lemma for Theorem 2.9 by a similar
idea as the one used in the proof of Lemma 2.5. By this lemma, if
an overpass of a knot diagram $D$ with respect to an
over-underpass sequence of $D$ crosses an adjacent underpass or an
underpass crosses an adjacent overpass, then the number of
crossings of $D$ is not minimal.

\begin{lm} If $D$ is a minimal diagram of a knot $K$ with respect to
crossings, $b(D)=k \geq 2$, and $s_1,f_1,s_2,f_2,...,s_k,f_k$ is
an over-underpass sequence of $D$, then no overpass of $D$ with
respect to $s_1,f_1,s_2,f_2,...,s_k,f_k$ crosses its adjacent
underpasses and no underpass crosses its adjacent overpasses,
where the adjacent underpasses of the $i$-th overpass $o_i$ are
the $(i-1)$-th and the $i$-th underpasses $u_{i-1}$ and $u_i$ for
each $i \in \{2,...,k\}$; the adjacent underpasses of $o_1$ are
$u_k$ and $u_1$; the adjacent overpasses of $u_i$ are $o_i$ and
$o_{i+1}$ for each $i \in \{1,...,k-1\}$; the adjacent overpasses
of $u_k$ are $o_k$ and $o_1$.
\end{lm}

\begin{proof} We may assume that each overpass of $D$ with respect to
$s_1,f_1,s_2,f_2,...,s_k,f_k$ crosses each underpass at most once
by Lemma 2.5.

Suppose that there is $i \in \{1,...,k\}$ such that the $i$-th
overpass $o_i=[s_i,f_i]$ of $D$ with respect to
$s_1,f_1,s_2,f_2,...,s_k,f_k$ crosses the $i$-th underpass
$u_i=[f_i,s_{i+1}]$, with the subscript counted modulo $k$. Let
$x$ be the crossing of $D$ between $o_i$ and $u_i$. Let $y$ be the
crossing of $D$ such that $y_+$ is the overcrossing of $o_i$ just
before $f_i$, and let $z$ be the crossing of $D$ such that $z_-$
is the undercrossing of $u_i$ just after $f_i$. Note that $x$,
$y$, $z$ need not be distinct. However, if some of them are
identical, then, by Reidemeister moves, we can reduce the crossing
easily. Let $r$ be the number of overcrossings on the arc
$\overrightarrow{x_+y_+}^{o_i}$ of $o_i$ from $x_+$ to $y_+$, and
let $s$ be the number of undercrossings on the arc
$\underrightarrow{z_-x_-}_{u_i}$ of $u_i$ from $z_-$ to $x_-$.
Temporarily, we regard the knot diagram $D$ as its regular
projection as we did in the the proof of Lemma 2.5. Let us take an
$\epsilon$-neighborhood $U_{o_i,\epsilon}$ of
$\overrightarrow{xy}^{o_i}$ and an $\epsilon$-neighborhood
$V_{u_i,\epsilon}$ of $\underrightarrow{zx}_{u_i}$ for a
sufficiently small positive real number $\epsilon$ as described in
the proof of Lemma 2.5. Let $a_1$ and $a_2$ be the first and the
second points from $s_i$ in $Bd(U_{o_i,\epsilon/2}) \cap o_i$,
respectively, and let $a_3$ and $a_4$ be the first and the second
points from $s_i$ in $Bd(U_{o_i,\epsilon/2}) \cap u_i$,
respectively. Similarly, let $b_1$ and $b_2$ be the first and the
second points from $s_i$ in $Bd(V_{u_i,\epsilon/2}) \cap o_i$,
respectively, and let $b_3$ and $b_4$ be the first and the second
points from $s_i$ in $Bd(V_{u_i,\epsilon/2}) \cap u_i$,
respectively. Let $\overline{a_2a_4}$ be the arc of
$Bd(U_{o_i,\epsilon/2})$ from $a_2$ to $a_4$ such that
$|\overline{a_2a_4} \cap D|=(r-1)+2$, and let $\overline{b_1b_3}$
be the arc of $Bd(V_{u_i,\epsilon/2})$ from $b_1$ to $b_3$ such
that $|\overline{b_1b_3} \cap D|=(s-1)+2$.

From now on, $D$ is the diagram of $K$ again. Let
$\underrightarrow{a_2a_4}$ be a regularly projecting arc in $R^3$
whose regular projection is $\overline{a_2a_4}$ such that
$\underrightarrow{a_2a_4}$ has no overcrossing, and let
$\overrightarrow{b_1b_3}$ be a regularly projecting arc in $R^3$
whose regular projection is $\overline{b_1b_3}$ such that
$\overrightarrow{b_1b_3}$ has no undercrossing. Then

Case 1. $r \leq s$: Let $D'=(D-(\overrightarrow{a_2f_i}^{o_i} \cup
\underrightarrow{f_ia_4}_{u_i})) \cup \underrightarrow{a_2a_4}$.
Then $D'$ is a diagram of $K$ and $c(D')=c(D)-s+(r-1)$, hence,
$c(D)=c(D')+(s-r)+1>c(D')$.

Case 2. $r>s$: Let $D'=(D-(\overrightarrow{b_1f_i}^{o_i} \cup
\underrightarrow{f_ib_3}_{u_i})) \cup \overrightarrow{b_1b_3}$.
Then $D'$ is a diagram of $K$ and $c(D')=c(D)-r+(s-1)$, hence,
$c(D)=c(D')+(r-s)+1>c(D')$.

Therefore, we have a diagram $D'$ of $K$ such that $c(D')<c(D)$.
This is a contradiction to $c(D)=c(K)$.

\bigskip
\centerline{\epsfxsize=3.5in\epsfbox{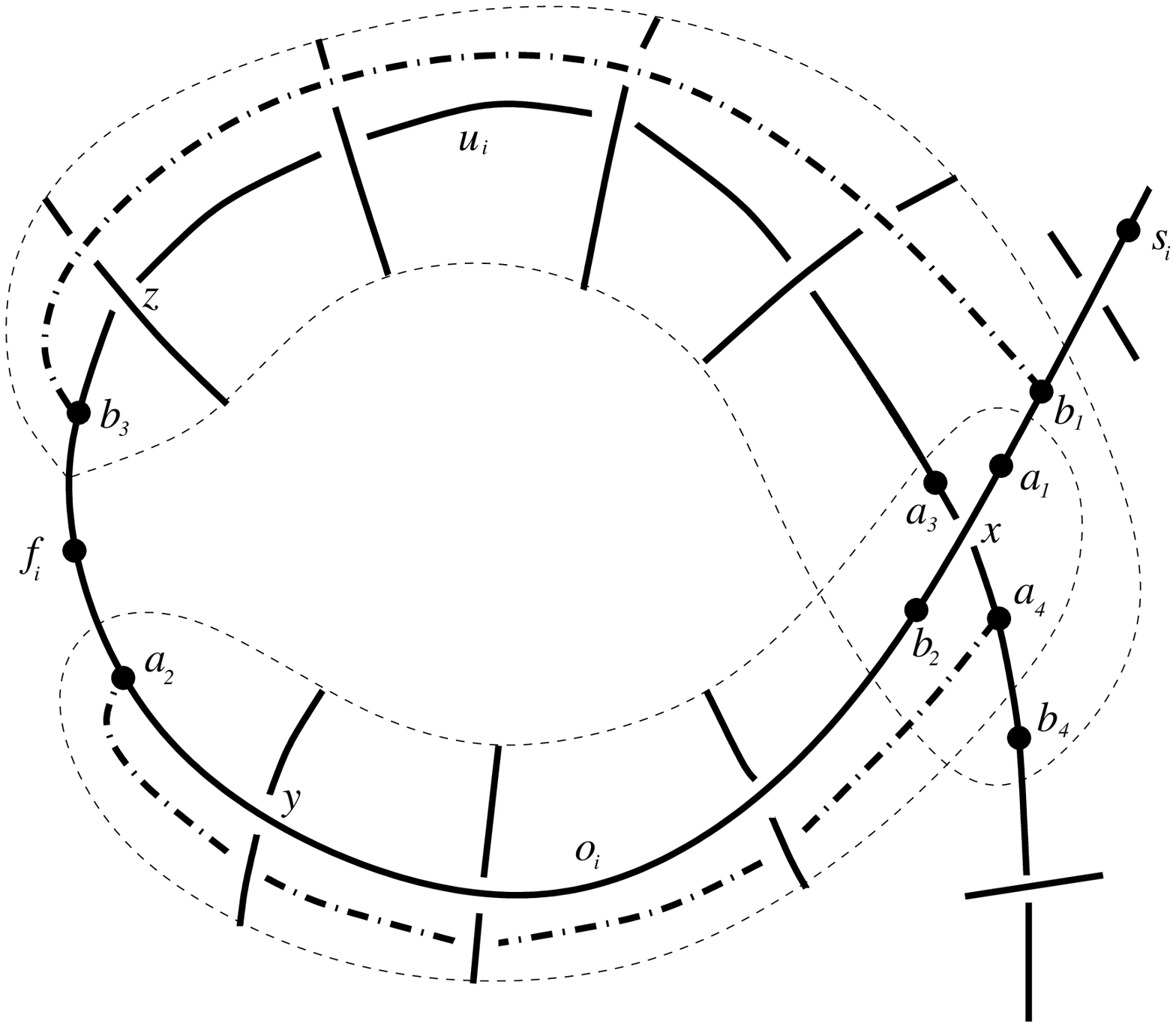}}
\medskip
\centerline{\small Figure 2. Proof of Lemma 2.6.}
\bigskip

As the other case, if there is $i \in \{1,...,k\}$ such that the
$i$-th underpass $u_i=[f_i,s_{i+1}]$ of $D$ with respect to
$s_1,f_1,s_2,f_2,...,s_k,f_k$ crosses the $(i+1)$-th overpass
$o_{i+1}=[s_{i+1},f_{i+1}]$, with the subscript counted modulo
$k$, then we can construct a diagram $D'$ of $K$ such that
$c(D')<c(D)$ by the same argument above. This proves the lemma.
\end{proof}

The following two lemmas are immediate consequences of Lemma 2.6
and the other keys for Theorem 2.9.

\begin{lm} If $D$ is a minimal diagram of a knot $K$ with respect to
crossings, $b(D)$ is a positive integer $k$, and
$s_1,f_1,s_2,f_2,...,s_k,f_k$ is an over-underpass sequence of
$D$, then there is no overpass of $D$ with respect to
$s_1,f_1,s_2,f_2,...,s_k,f_k$ which crosses all underpasses.
\end{lm}

\begin{proof} Suppose that an overpass of $D$ with respect to
$s_1,f_1,s_2,f_2,...,s_k,f_k$ crosses all underpasses. Then the
overpass crosses its adjacent underpasses. This is a contradiction
to $c(D)=c(K)$ by Lemma 2.6.
\end{proof}

\begin{lm} If $D$ is a minimal diagram of a knot $K$ with respect to
crossings, $b(D)=k \geq 2$, and $s_1,f_1,s_2,f_2,...,s_k,f_k$ is
an over-underpass sequence of $D$, then every overpass of $D$ with
respect to $s_1,f_1,s_2,f_2,...,s_k,f_k$ crosses at most $k-2$
underpasses.
\end{lm}

\begin{proof} Suppose that an overpass $o$ of $D$ with respect to
$s_1,f_1,s_2,f_2,...,s_k,f_k$ crosses $k-1$ underpasses. Then, by
Lemma 2.5, $o$ crosses each of the $k-1$ underpasses exactly once.
Since every overpass has two adjacent underpasses and $o$ crosses
$k-1$ underpasses, $o$ must cross at least one of its adjacent
underpasses. This is a contradiction to $c(D)=c(K)$ by Lemma 2.6.
\end{proof}

\begin{thm} If $D$ is a minimal diagram of a knot $K$ with respect to
crossings, then $b(D) \leq c(D) \leq b(D)(b(D)-2)$.
\end{thm}

\begin{proof} Let $D$ be a diagram of a knot $K$ such that
$c(D)=c(K)$. If $K$ is a trivial knot, then $D$ is a simple closed
curve on the plane, hence, $c(D)=b(D)=0$. Remark that $b(D)$ can
not be $1$. If $b(D)=1$, then, by Corollary 2.2, $D$ represents a
trivial knot, i.e., $K$ is a trivial knot, hence, $c(D)=c(K)=0$.
This contradicts to Lemma 2.4. Next, we claim that $b(D)$ can not
be 2 if $D$ is a minimal diagram of a nontrivial knot $K$ with
respect to crossings. This follows from the same argument in the
proof of Lemma 2.8. Suppose that $K$ is a nontrivial knot,
$c(D)=c(K)$, $b(D)=2$, and $s_1,f_1,s_2,f_2$ is an over-underpass
sequence of $D$. Then each overpass of $D$ with respect to
$s_1,f_1,s_2,f_2$ must cross one of the underpasses. However, all
underpasses are adjacent to each overpass. Suppose that $b(D)=k
\geq 3$ and $s_1,f_1,s_2,f_2,...,s_k,f_k$ is an over-underpass
sequence of $D$. We may assume that the overpasses
$[s_1,f_1],[s_2,f_2],...,[s_{k-1}, f_{k-1}],[s_k,f_k]$ are
disjoint closed arcs on the projection plane of $D$ and the
underpasses $[f_1,s_2],[f_2,s_3],...,[f_{k-1},s_k],[f_k,s_1]$ are
also disjoint closed arcs on it. Remark that the number of
crossings of $D$ is the number of intersections of projections of
overpasses and underpasses of $D$. Since we have $k$ overpasses
and $k$ underpasses, by Lemma 2.5 and Lemma 2.8, there are at most
$k(k-2)$ crossings among overpasses and underpasses.
\end{proof}

Notice that $c(K)\leq b(K)(b(K)-2)$ does not hold in general. An
example can be given with $K$ being the $(2,p)$-torus knot. We
have $b(K)=2$ in this case.

The following theorem is an immediate consequence of Theorem 2.9
if $K$ is nontrivial. It means that the number of overpasses of a
minimal knot diagram with respect to crossings can also be
estimated by the number of crossings. Once more, in this case, the
number of overpasses needs not be minimal and is at least 3 as
shown in the proof of Theorem 2.9.

\begin{thm} If $K$ is a nontrivial knot and $D$ is a minimal diagram
of $K$ with respect to crossings, then $1+ \sqrt{1+c(D)} \leq b(D)
\leq c(D)$.
\end{thm}

\section{The most nonalternating knots are $(k-1,\pm k)$-torus knots}

Theorem 2.10 provides us the most optimal lower bound for the
number of overpasses of $D$ when $D$ is a minimal knot diagram
with respect to crossings. Here, an interesting problem occurs. In
Lemma 2.4, we showed that a knot diagram $D$ is alternating if and
only if $c(D)=b(D)$. Also, by Kauffman \cite{K}, an alternating
knot has a minimal knot diagram with respect to crossings which is
alternating. Hence, the second inequality of Theorem 2.10 becomes
an equality if and only if the knot is an alternating knot. On the
other hand, as shown by Murasugi \cite{MUR}, the standard diagram
$D$ of a $(k-1,\pm k)$-torus knot is a minimal knot diagram with
respect to crossings and $c(D)=k(k-2)$. So the first inequality of
Theorem 2.10 becomes an equality for the $(k-1,\pm k)$-torus knot
with $b(D)=k$. Is the converse true under the condition that $D$
is a minimal knot diagram with respect to crossings? We prove this
as the last theorem of this paper.

\begin{thm} If $D$ is a minimal diagram of a nontrivial knot $K$
with respect to crossings and $c(D)=b(D)(b(D)-2)$, then $D$ is the
standard diagram of either the $(b(D)-1,b(D))$-torus knot or the
$(b(D)-1,-b(D))$-torus knot. Hence, $K$ is the $(b(D)-1,\pm
b(D))$-torus knot.
\end{thm}

\begin{proof} Let $D$ be a minimal diagram of a nontrivial knot $K$
with respect to crossings such that $c(D)=b(D)(b(D)-2)$. Suppose
that $b(D)=k$. Then $k \geq 3$ (See the proof of Theorem 2.9). In
order to prove this theorem, we will consider all possible knot
diagrams satisfying the hypothesis and show that they can only be
the standard diagrams of either the $(k-1,k)$-torus knot or the
$(k-1,-k)$-torus knot depending on signs of crossings. For
convenience, we regard the knot diagram $D$ as its regular
projection here. Also, to visualize overpasses and underpasses
clearly, we imagine blue and red colors for overpasses and
underpasses, respectively. In this sense, crossings are only the
intersections of blue arcs and red arcs on the plane, i.e.,
crossings are purple!

Now, let us start drawing all possible knot diagrams satisfying
the hypothesis. Notice that such knot diagrams must satisfy the
following 3 rules:

Rule 1: Every overpass and every underpass intersect each
underpass and each overpass at most once, respectively (By Lemma
2.5);

Rule 2: No overpass and no underpass intersect its adjacent
underpasses and its adjacent overpasses, respectively (By Lemma
2.6);

Rule 3: Every overpass and every underpass intersect $k-2$
underpasses and $k-2$ overpasses, respectively (By Lemma 2.8).

If some overpass has less than $k-2$ overcrossings, then at least
one overpass has more than $k-2$ overcrossings; similarly, if some
underpass has less than $k-2$ undercrossings, then at least one
underpass has more than $k-2$ undercrossings. This is a
contradiction to Lemma 2.8. From now on, we describe drawing the
knot diagrams. Let us take a point $s_1$ on the plane and draw the
first overpass $o_1=[s_1,f_1]$. Draw the first underpass
$u_1=[f_1,s_2]$ and the second overpass $o_2=[s_2,f_2]$. Until
now, we can not have any crossing by Rule 2. When we draw the
second underpass $u_2$ from $f_2$, we must make $u_2$ intersect
$o_1$ first by Rule 2,3. Notice that we have only 2 cases that
$u_2$ intersects $o_1$ first as follows.

Case 1. The sign $sign(u_2,o_1)$ of crossing between $u_2$ and
$o_1$ is $-1$,

where $sign(u_2,o_1)=+1=sign(o_1,u_2)$ if $o_1$ intersects $u_2$
from left to right and $sign(u_2,o_1)=-1=sign(o_1,u_2)$ if $o_1$
intersects $u_2$ from right to left when we look $u_2$ as a line
segment passing through the crossing upward. Notice that, by Rule
1, we can define the sign of each crossing by this way.

After $u_2$ intersects $o_1$ first with $sign(u_2,o_1)=-1$, we
must change color from red to blue by Rule 1,2. Take a point $s_3$
so that $u_2=[f_2,s_3]$ has only one undercrossing. Then the third
overpass $o_3$ must intersect $u_1$ first with $sign(o_3,u_1)=-1$
by Rule 2,3. After $o_3$ intersects $u_1$ first, we must change
color from blue to red by Rule 1,2. Take a point $f_3$ so that
$o_3=[s_3,f_3]$ has only one overcrossing. Then the third
underpass $u_3$ can intersect only either $o_1$ first with
$sign(u_3,o_1)=-1$ or $o_2$ first with $sign(u_3,o_2)=-1$. We
claim that $u_3$ must intersect $o_2$ first with
$sign(u_3,o_2)=-1$. Suppose that $u_3$ intersects $o_1$ first with
$sign(o_3,u_1)=-1$. Then, after $u_3$ intersects $o_1$ first, we
must change color from red to blue by Rule 1,2. Take a point $s_4$
so that $u_3=[f_3,s_4]$ has only one undercrossing. In this case,
we can not draw any knot diagram such that $u_3$ intersects $o_2$.
This is a contradiction to Rule 2,3. Hence, $u_3$ intersects $o_2$
first with $sign(u_3,o_2)=-1$. After $u_3$ intersects $o_2$ first,
take a point $x_3$ which is neither $s_1$ nor a crossing so that
the arc $[f_3,x_3]$ of $u_3$ has only one undercrossing. Notice
that we can connect $x_3$ and $s_1$ by an arc without crossing.

\bigskip
\centerline{\epsfxsize=4.8in\epsfbox{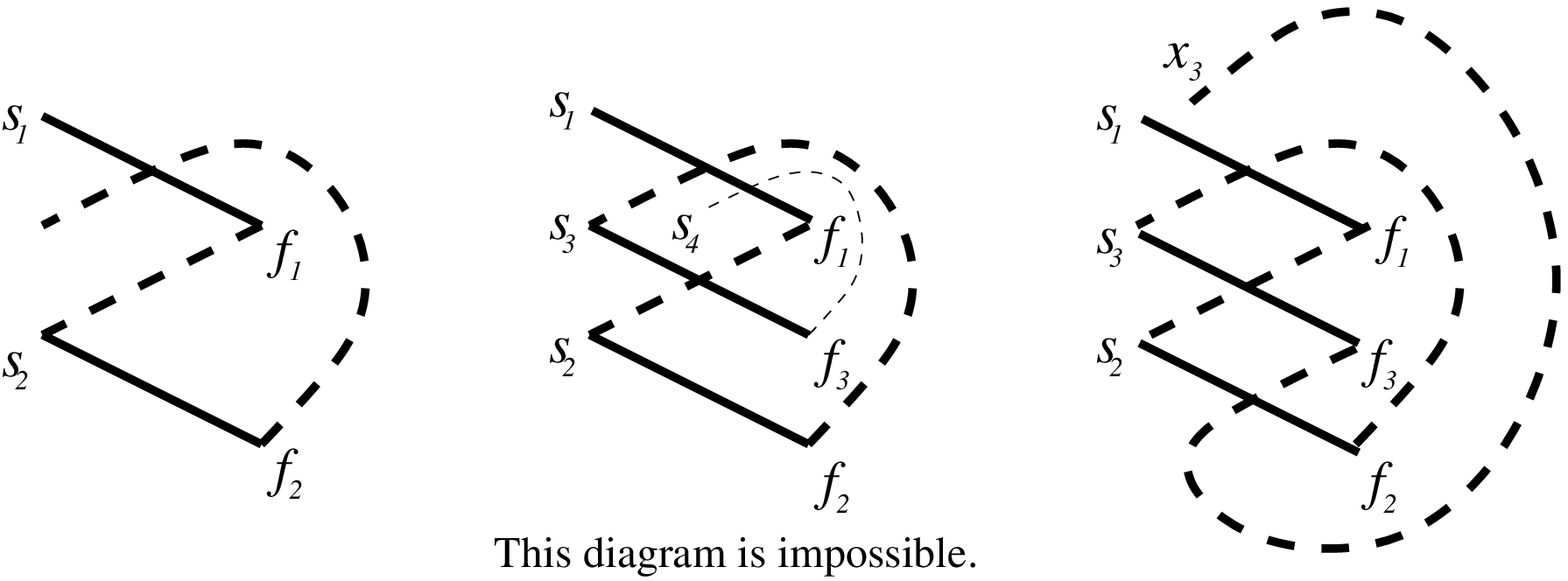}}
\medskip
\centerline{\small Figure 3. The first three overpasses and
underpasses, $sign(u_2,o_1)=-1$.}
\bigskip

If we have only 3 overpasses and 3 underpasses, we must connect
$x_3$ and $s_1$ by an arc without crossing. This is the only way
to draw the knot diagrams satisfying Rule 1,2,3 and gives us the
standard diagram of the $(2,-3)$-torus knot.

See Figure 3, where solid arcs are supposed to be of the blue
color (overpasses) and dashed arcs of the red color (underpasses).

Suppose that we have more than 3 overpasses, i.e., $k>3$. Then the
third underpass $u_3$ must intersect $o_1$ second by Rule 1,2,3
and we have only 2 cases that $u_3$ intersects $o_1$ second.

\noindent{\bf Claim:} $u_3$ must intersect $o_1$ second with
$sign(u_3,o_1)=-1$. See Figure 4(i).

\noindent{\it Proof of the Claim:} Suppose that $u_3$ intersects
$o_1$ second with $sign(u_3,o_1)=+1$. See Figure 4(ii). Then,
after $u_3$ intersects $o_1$ second, we must change color from red
to blue by Rule 1. Take a point $s_4$ so that $u_3=[f_3,s_4]$ has
only two undercrossings. Then $o_4$ must intersect $u_2$ first
with $sign(o_4,u_2)=-1$ and $u_1$ second with $sign(o_4,u_1)=+1$
by Rule 1,2,3. See Figure 4(iii).

After $o_4$ intersects $u_1$ second, we must change color from
blue to red by Rule 1. Take a point $f_4$ so that $o_4=[s_4,f_4]$
has only two overcrossings. Then $u_4$ can intersect only either
$o_1$ first with $sign(u_4,o_1)=+1$ or $o_3$ first with
$sign(u_4,o_3)=-1$. However, $u_4$ can not intersect $o_1$ first.
If $u_4$ intersects $o_1$ first with $sign(u_4,o_1)=+1$, we must
change color from red to blue by Rule 1,2. By a similar argument
as the one used before, in this case, we can not draw any knot
diagram such that $u_4$ intersects $o_2$. This is a contradiction
to Rule 2,3. Hence, $u_4$ must intersect $o_3$ first with
$sign(u_4,o_3)=-1$.

After $u_4$ intersects $o_3$ first, $u_4$ can intersect only
either $o_1$ second with $sign(u_4,o_1)=+1$ or $o_2$ second with
$sign(u_4,o_2)=+1$. However, $u_4$ can not intersect $o_1$ second.
If $u_4$ intersects $o_1$ second with $sign(u_4,o_1)=+1$, after
$u_4$ intersects $o_1$ second, we must change color from red to
blue by Rule 1,2. In this case, we can not draw any knot diagram
such that $u_4$ intersects $o_2$. This is a contradiction to Rule
2,3. Hence, the only way to draw $u_4$ is that $u_4$ intersects
$o_3$ first with $sign(u_4,o_3)=-1$ and $o_2$ second with
$sign(u_4,o_2)=+1$. See Figure 4(iv).

\bigskip
\centerline{\epsfxsize=4.5in\epsfbox{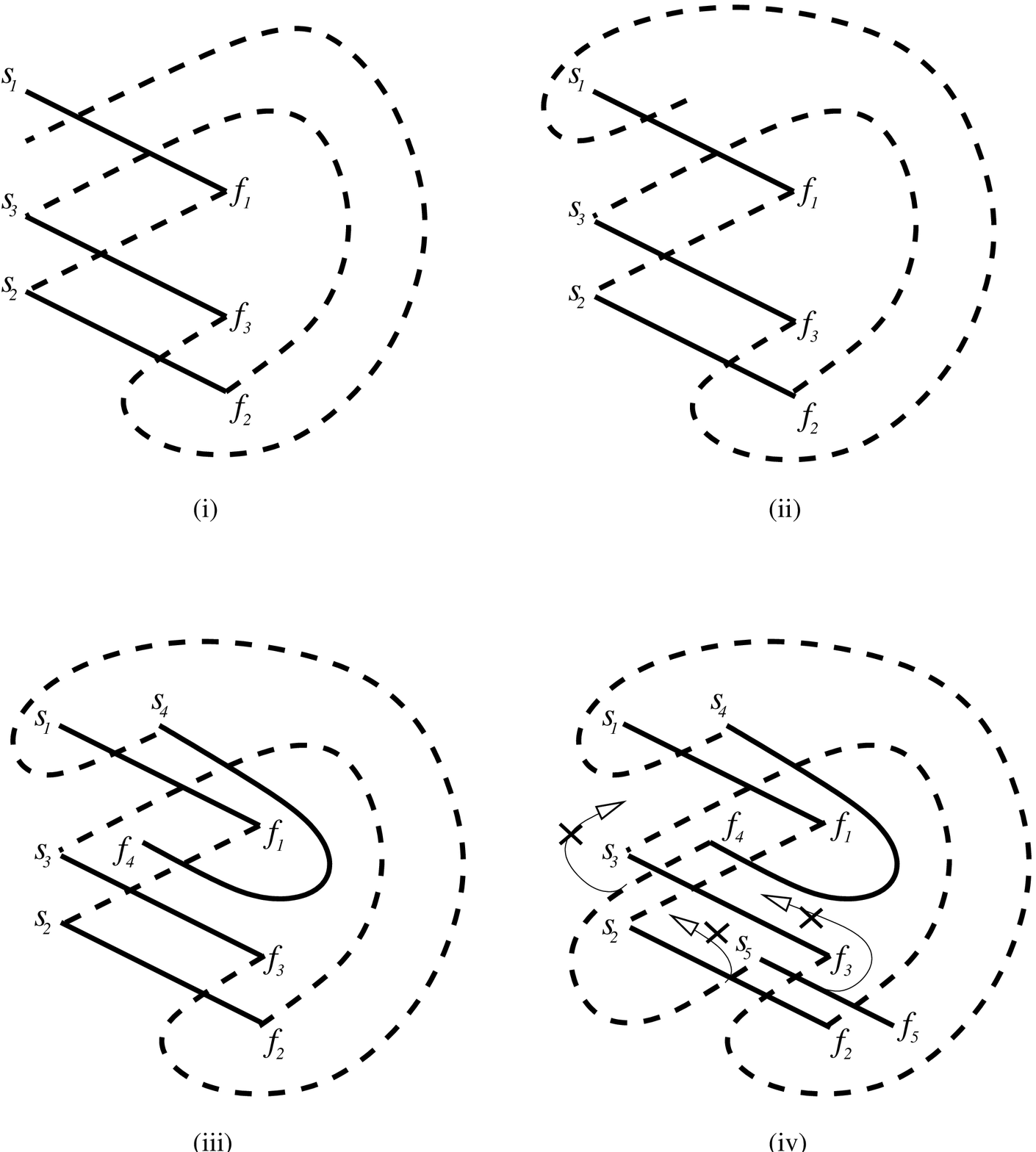}}
\medskip
\centerline{\small Figure 4. The first three overpasses and
underpasses can only be drawn as in (i).}
\bigskip

After $u_4$ intersects $o_2$ second, we must change color from red
to blue by Rule 1. Take a point $s_5$ so that $u_4=[f_4,s_5]$ has
only two undercrossings. Then $o_5$ can intersect only either
$u_1$ first with $sign(o_5,u_1)=+1$ or $u_3$ first with
$sign(o_5,u_3)=-1$. However, $o_5$ can not intersect $u_1$ first.
If $o_5$ intersects $u_1$ first with $sign(o_5,u_1)=+1$, then we
must change color from blue to red by Rule 1,2. In this case, we
can not draw any knot diagram such that $o_5$ intersects $u_2$.
This is a contradiction to Rule 2,3. Hence, $o_5$ must intersect
$u_3$ first with $sign(o_5,u_3)=-1$. Again, see Figure 4(iv).

After $o_5$ intersects $u_3$ first, $o_5$ can intersect only
either $u_1$ second with $sign(o_5,u_1)=+1$ or $u_2$ second with
$sign(o_5,u_2)=+1$. However, $o_5$ can not intersect $u_1$ second.
If $o_5$ intersects $u_1$ second with $sign(o_5,u_1)=+1$, after
$o_5$ intersects $u_1$ second, we must change color from blue to
red by Rule 1,2. In this case, we can not draw any knot diagram
such that $o_5$ intersects $u_2$. This is a contradiction to Rule
2,3. Hence, $o_5$ must intersect $u_2$ second with
$sign(o_5,u_2)=+1$. Once more, see Figure 4(iv).

After $o_5$ intersects $u_2$ second, we must change color from
blue to red by Rule 1. In this case, we can not draw any knot
diagram such that $o_5$ intersects $u_1$. This is a contradiction
to Rule 2,3. Hence, when $sign(u_3,o_1)=+1$, we can not draw any
knot diagram satisfying Rule 1,2,3. Therefore, $sign(u_3,o_1)$
must be $-1$. This finishes the proof of the Claim.

Now, we use an induction on the order of the over-underpass
sequence. The argument used here is just a generalization of the
previous one. Hence, our argument will be slightly sketchy here.

Suppose, without loss of generality, that $3 \leq n<k$ and
$s_1,f_1,s_2,f_2,...,s_n,f_n$ is an over-underpass sequence of the
standard diagram of the $(n-1,-n)$-torus knot.

\bigskip
\centerline{\epsfysize=5.5in\epsfbox{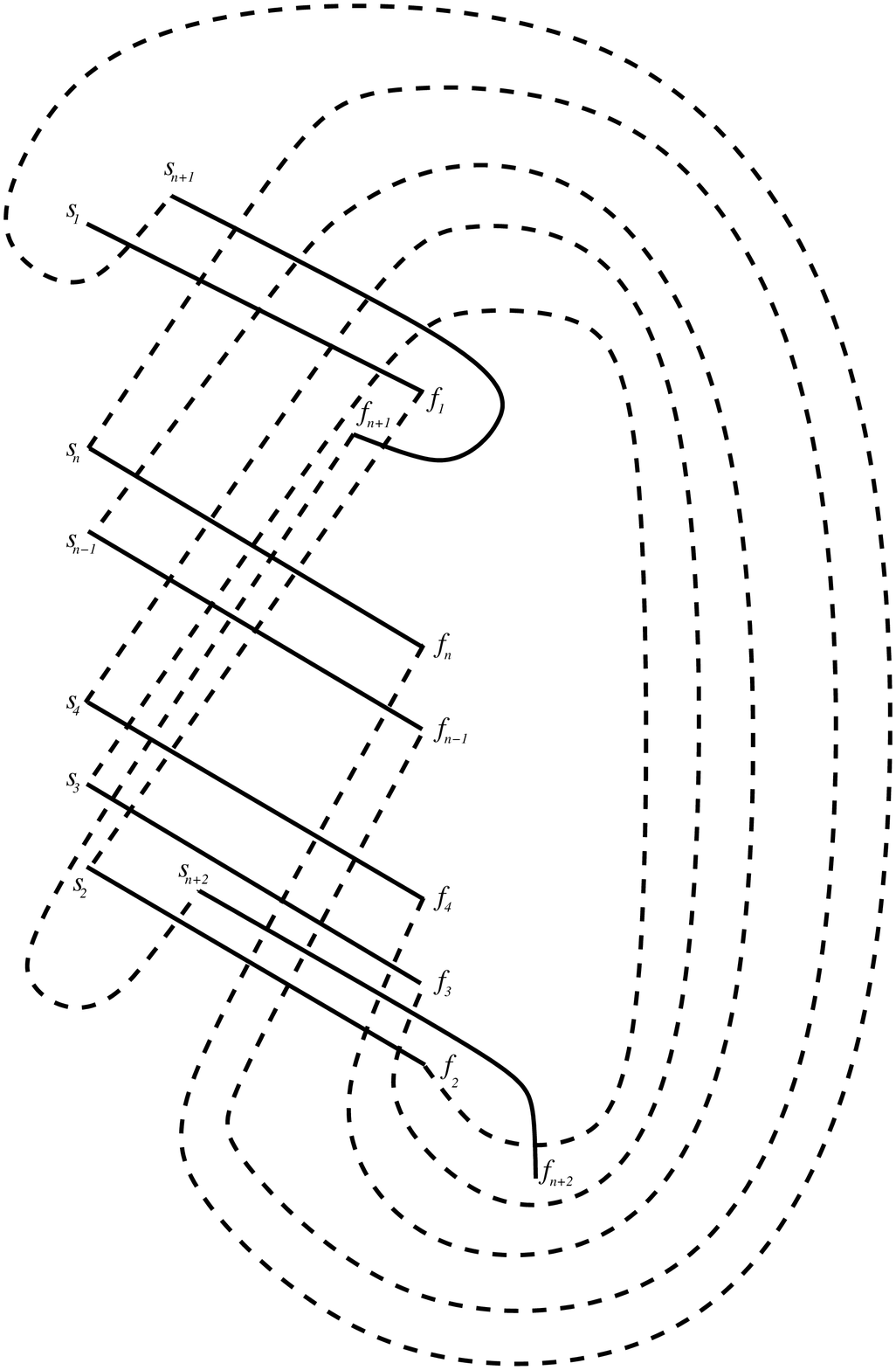}}
\medskip
\centerline{\small Figure 5. It is impossible for $u_n$ to
intersect $o_1$ with $sign(u_n,o_1)=-1$.}
\bigskip

Let $o_i=[s_i,f_i]$ and $u_i=[f_i,s_{i+1}]$ for each $i \in
\{1,...,n-1\}$, and let $o_n=[s_n,f_n]$ and $u^*_n=[f_n,s_1]$.
Take a point $x_n$ of $u^*_n$ between the $(n-2)$-th undercrossing
of $u^*_n$ and $s_1$. Then the arc $[f_n,x_n]$ of $u^*_n$ has also
$n-2$ undercrossings.

Notice that $o_1,u_1,o_2,u_2,...,o_{n-1},u_{n-1},o_n$ are
overpasses and underpasses of a knot diagram $D$ satisfying Rule
1,2,3 and $[f_n,x_n]$ is an arc of the $n$-th underpass $u_n$ of
$D$. Suppose that $D$ has more than $n$ overpasses. Then there is
an $((n-2)+1)$-th undercrossing on $u_n$ by $o_1$ according to
Rule 1,2,3. We claim that at this undercrossing, $u_n$ must
intersect $o_1$ with $sign(u_n,o_1)=-1$.

To prove this claim, suppose that $u_n$ intersects $o_1$  with
$sign(u_n,o_1)=+1$. See Figure 5. Then, after $u_n$ intersects
$o_1$ with $sign(u_n,o_1)=+1$, we must change color from red to
blue by Rule 1. Take a point $s_{n+1}$ so that $u_n=[f_n,s_{n+1}]$
has $(n-2)+1$ undercrossings. Then the only way to draw $o_{n+1}$
is that $o_{n+1}$ intersects $u_{n-1}$ first, $u_{n-2}$ second,
..., $u_2$ $(n-2)$-th, and $u_1$ $((n-2)+1)$-th so that
$sign(o_{n+1},u_{n-1})=sign(o_{n+1},u_{n-2})=...=sign(o_{n+1},u_2)=-1$,
and $sign(o_{n+1},u_1)=+1$ by Rule 1,2,3. After $o_{n+1}$
intersects $u_1$ with $sign(o_{n+1},u_1)=+1$, we must change color
from blue to red by Rule 1. Take a point $f_{n+1}$ so that
$o_{n+1}$ has $(n-2)+1$ overcrossings. Then the only way to draw
$u_{n+1}$ is that $u_{n+1}$ intersects $o_n$ first, $o_{n-1}$
second, ..., $o_3$ $(n-2)$-th, and $o_2$ $((n-2)+1)$-th so that
$sign(u_{n+1},o_n)=sign(u_{n+1},o_{n-1})=...=sign(u_{n+1},o_3)=-1$,
and $sign(u_{n+1},o_2)=+1$ by Rule 1,2,3. After $u_{n+1}$
intersects $o_2$ with $sign(u_{n+1},o_2)=+1$, we must change color
from red to blue by Rule 1. Take a point $s_{n+2}$ so that
$u_{n+1}$ has $(n-2)+1$ undercrossings. Then the only way to draw
$o_{n+2}$ is that $o_{n+2}$ intersects $u_n$ first, $u_{n-1}$
second, ..., $u_3$ $(n-2)$-th, and $u_2$ $((n-2)+1)$-th so that
$sign(o_{n+2},u_n)=sign(o_{n+2},u_{n-1})=...=sign(o_{n+2},u_3)=-1$,
and $sign(o_{n+2},u_2)=+1$ by Rule 1,2,3. After $o_{n+2}$
intersects $u_2$ with $sign(o_{n+2},u_2)=+1$, we must change color
from blue to red by Rule 1. However, in this case, we can not draw
any knot diagram such that $o_{n+2}$ intersects $u_1$. This is a
contradiction to Rule 2,3. Therefore, $sign(u_n,o_1)$ can not be
$+1$.

As a next step, we claim that the only way to draw $o_{n+1}$ and
$u_{n+1}$ gives us the standard diagram of the
$((n+1)-1,-(n+1))$-torus knot (see Figure 6). This will finish the
induction.

After $u_n$ intersects $o_1$ $((n-2)+1)$-th with
$sign(u_n,o_1)=-1$, we must change color from red to blue by Rule
1,2. Take a point $s_{n+1}$ so that $u_n=[f_n,s_{n+1}]$ has only
$(n-2)+1$ undercrossings. Then the only way to draw $o_{n+1}$ is
that $o_{n+1}$ intersects $u_{n-1}$ first, $u_{n-2}$ second, ...,
$u_1$ $((n-2)+1)$-th so that
$sign(o_{n+1},u_{n-1})=sign(o_{n+1},u_{n-2})=...=sign(o_{n+1},u_1)=-1$
by Rule 1,2,3. Suppose that $o_{n+1}$ intersects $u_i$ first for
some $i<n-1$. Then $sign(o_{n+1},u_i)$ must be $-1$ and $o_{n+1}$
must intersect $u_i$ first, $u_{i-1}$ second, ..., $u_1$ $i$-th so
that
$sign(o_{n+1},u_i)=sign(o_{n+1},u_{i-1})=...=sign(o_{n+1},u_1)=-1$
by Rule 1,2,3. However, after $o_{n+1}$ intersects $u_1$ $i$-th
with $sign(o_{n+1},u_1)=-1$, $o_{n+1}$ must stop before
intersecting $u_n$ and we must change color from blue to red by
Rule 1,2. In this case, we can not draw any knot diagram such that
$o_{n+1}$ intersects $u_{n-1}$. This is a contradiction to Rule
2,3. Hence, $o_{n+1}$ must intersect $u_{n-1}$ first, $u_{n-2}$
second, ..., $u_1$ $((n-2)+1)$-th so that
$sign(o_{n+1},u_{n-1})=sign(o_{n+1},u_{n-2})=...=sign(o_{n+1},u_1)=-1$.
Take a point $f_{n+1}$ so that $o_{n+1}=[s_{n+1},f_{n+1}]$ has
only $(n-2)+1$ overcrossings. Then, by a similar argument as
before, we can show that the only way to draw $u_{n+1}$ is that
$u_{n+1}$ intersects $o_n$ first, $o_{n-1}$ second, ..., $o_2$
$((n-2)+1)$-th so that
$sign(u_{n+1},o_n)=sign(u_{n+1},o_{n-1})=...=sign(u_{n+1},o_2)=-1$
by Rule 1,2,3.

\bigskip
\centerline{\epsfysize=6in\epsfbox{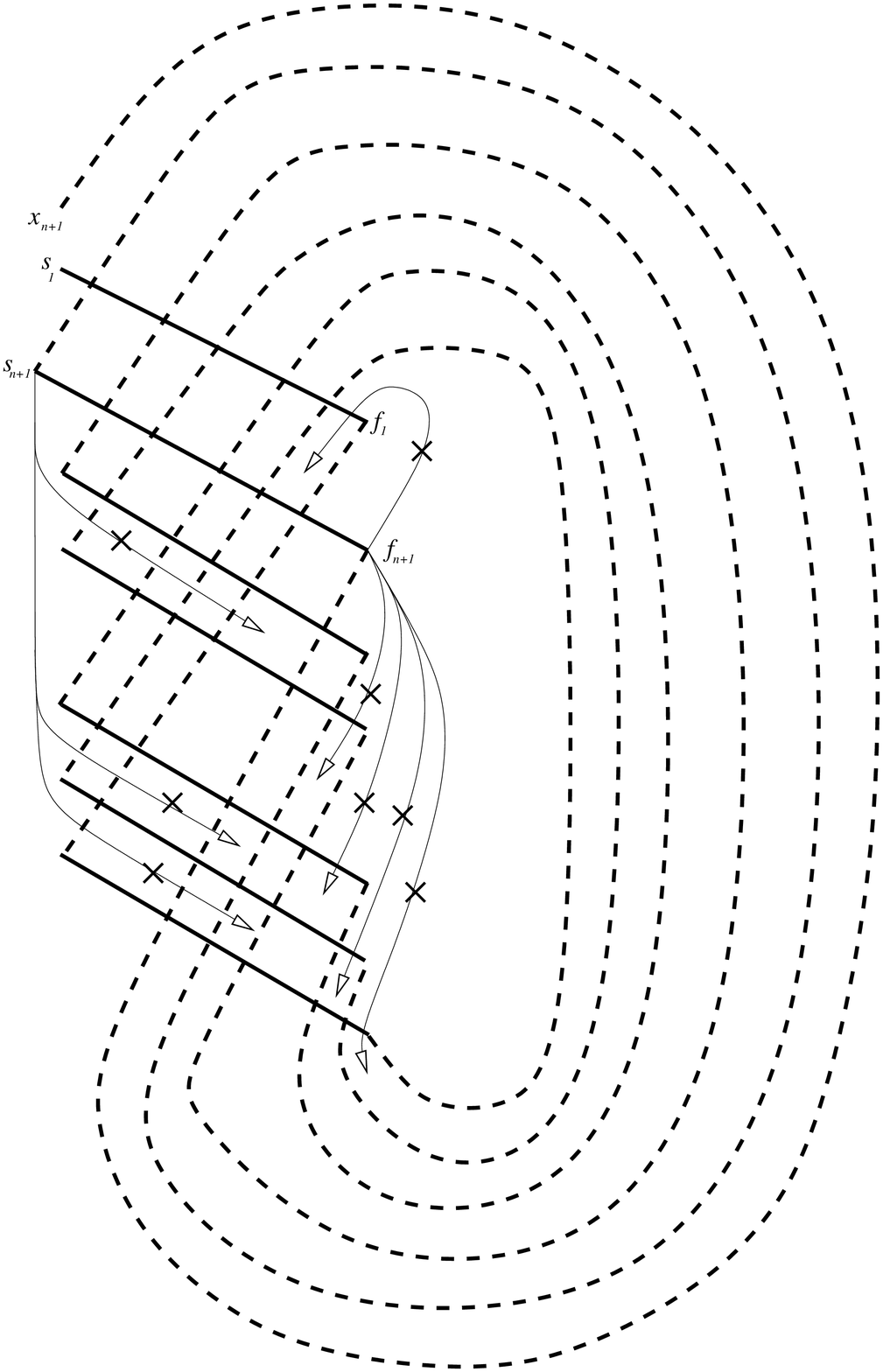}}
\medskip
\centerline{\small Figure 6. The only way to draw $o_{n+1}$ and
$u_{n+1}$.}
\bigskip

Suppose that $u_{n+1}$ intersects $o_i$ first for some $i<n$. Then
$sign(u_{n+1},o_i)$ must be $-1$ and $u_{n+1}$ must intersect
$o_i$ first, $o_{i-1}$ second, ..., $o_1$ $i$-th so that
$sign(u_{n+1},o_i)=sign(u_{n+1},o_{i-1})=...=sign(u_{n+1},o_1)=-1$
by Rule 1,2,3. However, after $u_{n+1}$ intersects $o_1$ $i$-th
with $sign(u_{n+1},o_1)=-1$, $u_{n+1}$ must stop before
intersecting $o_{n+1}$ and we must change color from red to blue
by Rule 1,2. In this case, we can not draw any knot diagram such
that $u_{n+1}$ intersects $o_n$. This is a contradiction to Rule
2,3. Hence, $u_{n+1}$ must intersect $o_n$ first, $o_{n-1}$
second, ..., $o_2$ $((n-2)+1)$-th so that
$sign(u_{n+1},o_n)=sign(u_{n+1},o_{n-1})=...=sign(u_{n+1},o_2)=-1$.
Take a point $x_{n+1}$ so that the arc $[f_{n+1},x_{n+1}]$ of
$u_{n+1}$ has only $(n-2)+1$ undercrossings. Notice that, in this
case, we can connect $x_{n+1}$ and $s_1$ by an arc without
crossing so that we complete drawing the knot diagram.

This is the only way to draw a knot diagram with $n+1$ overpasses
satisfying the Rules 1,2,3.

To complete the argument, when $n=k-1$, we can get the standard
diagram $D$ of the $(k-1,-k)$-torus knot and this is the only way
to draw the knot diagrams satisfying Rule 1,2,3, and hence, the
theorem is proved in Case 1.

Case 2. The sign $sign(u_2,o_1)$ of crossing between $u_2$ and
$o_1$ is $+1$.

By the same argument as used in Case 1, we can get only standard
diagrams of the $(b(D)-1,b(D))$-torus knot. This proves the
theorem.
\end{proof}

\bigskip
\centerline{\epsfysize=3in\epsfbox{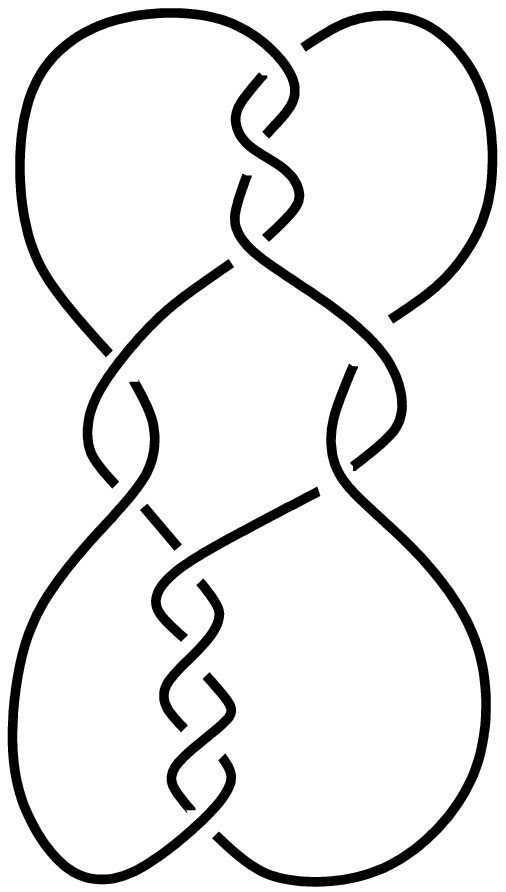}}
\medskip
\centerline{\small Figure 7. One of Goeritz's unknot diagrams.}
\bigskip

\section{A final remark}

One may view Lemmas 2.5 and 2.6 from another perspective. Suppose
$D$ is a knot diagram, not necessarily being minimal with respect
to crossings. If either an overpass interesects an underpass more
than once or an overpass interesects its adjacent underpass, then
Lemmas 2.5 and 2.6 give a specific way to change the knot diagram
$D$ by isotopy with the number of crossings reduced.

Consider a knot diagram $D$ satisfying the following two
conditions:

(1) each overpass intersects each underpass at most once; and

(2) each overpass does not intersect its adjacent underpasses.

For such a knot diagram $D$, we have
$$b(D)\leq c(D)\leq b(D)(b(D)-2).$$
When $c(D)=b(D)$, $D$ is an alternating knot diagram. By Theorem
3.1, when $c(D)=b(D)(b(D)-2)$, $D$ is the standard diagram of a
$(b(D)-1,\pm b(D))$ torus knot. Thus, it is natural to wonder
whether conditions (1) and (2) above are sufficient for a knot
diagram of a prime knot being minimal with respect to crossings.
Goeritz's unknot diagrams show that the answer to this question is
negative. See Figure 7. Are there any other necessary conditions
on the overpasses and underpasses of a knot diagram for it being
minimal with respect to crossings? This will be the topic of our
further investigation.

\end{document}